\documentclass{article}

\usepackage{amsmath,amssymb,amsthm}
\usepackage[margin=1in]{geometry} 
\usepackage{hyperref} 
\usepackage{graphicx} 
\usepackage{listings}

\title{The use of the invariant's properties in the\\ primality test and prime search.}
\author{Juan Hernández-Toro }
\date{\today}

\begin{document}

\maketitle

\begin{abstract}
The purpose of this article is to delve into the properties of invariants. The properties, explained in [2], reveal new ways to develop algorithms that allow us to test the primality of a number. In this article, some of these are shown, indicating the advantages and disadvantages of these new algorithms. The information provided by these algorithms also gives additional information regarding the factorization of a compound number.\\ 
\textbf{keywords:}Prime Numbers, Compound Numbers,  Primality test, Factorization.
\end{abstract}

\section{Introduction}
The concept of invariant could be used to explain the behavior of compound numbers. Using the invariant, it is possible to infer some additional properties from the compound numbers, such as non-trivial quadratic remainders and non-trivial triangular remainders.  If one of these properties is found, then it is automatically assumed that this number is not prime. Additionally, these remainders provide information regarding his factorization. In the conclusion of this article, some algorithms are shown to test primality using non-trivial invariants, non-trivial quadratic remainders, and non-trivial triangular remainders.
\section{Some definitions regarding the quadratic number's remainder}

\subsection{Definition of quadratic symmetric numbers  }\label{Cdef1}
The numbers whose quadratic powers have equivalent remainder are called quadratic symmetric numbers. For $p<m$, $s<m$ and $p\neq s$, $p$ is the Symmetric number of $s$ regarding $m$ if $s^2 \equiv p^2{\pmod {m}}$. This always happen if $m-s=p$.
\subsubsection{Demonstration:}
  
\begin{equation}
    p^{2}=m^{2}-m\cdot s+s^{2} \equiv s^{2}{\pmod {m}}
    \label{Ceq1}
\end{equation}
\subsubsection{Example:}
For $s=4$ and $m=15$ the symmetric number of n is $15-4=11$. As result $ 4^2 \equiv 11^2 \equiv 1{\pmod {15}}$. 

\subsection{Definition of numbers with trivial quadratic remainder}\label{Cdef2}
A number q with trivial quadratic remainder regarding  $m$ is defined as $q^2\equiv t^2 {\pmod {m}}$ and $q=t$.
\subsubsection{Conclusion:}
All the numbers $m > 1$ have Trivial quadratic remainder numbers. If $q^2<m$ then q is a trivial quadratic remainder. 
\subsubsection{Example:}
The numbers with trivial quadratic remainder regarding 15 are (0,1,2,3). 

\subsection{Definition of symmetric numbers with trivial quadratic remainder}\label{Cdef3}
If $t$ is a number with quadratic remainder regarding $m$ then $s$ is a Symmetric number with trivial quadratic remainder if $s=m-t$ 
\subsubsection{Conclusion:}
There are the same quantity of symmetric numbers with trivial quadratic remainder than numbers with trivial quadratic remainder.
\subsubsection{Example:}
The numbers with symmetric trivial quadratic remainder regarding 15 are (15,14,13,12).

\subsection{Definition of numbers with non trivial quadratic remainder}\label{Cdef4}
A number q with non trivial quadratic number remainder regarding $m$ is defined as $q^2\equiv d^2 {\pmod {m}}$ where $d^2<m$, $q\neq d$ and $m-q\neq d$

\subsubsection{Example:}
The numbers with non trivial quadratic regarding 15 are (4,7,8,11). 
 
\section{Some definitions regarding the triangular number's remainders}
\subsection{Definition of symmetric triangular numbers}\label{Tdef1}
The odd numbers whose triangular numbers have equivalent remainder are called triangular symmetric numbers. For $p<m$, $s<m$, m odd number and $p\neq s$; $p$ is the symmetric triangular number of $s$ regarding $m$ if $T(s)\equiv T(p){\pmod {m}}$.
\subsubsection{Demonstration:}
This always happen if $m-s-1=p$.  \footnote{For the even number the solution is slightly different }
\begin{equation}
\begin{split}
    \frac{p^2+p}{2}=\frac{m^2-2ms-2m+s^2+2s+1+m-s-1}{2}\\ \equiv \frac{s^2+2s+1-s-1}{2}=\frac{s^2+s}{2} {\pmod {m}}
    \label{Teq1}
		\end{split}
\end{equation}

\subsubsection{Example:}
For $s=6$ and $m=15$ the symmetric number of n is $15-6-1=8$. As result $ T(6)=21\equiv T(8)=36\equiv 6{\pmod {15}}$. 

\subsection{Definition of numbers with trivial triangular remainder}\label{Tdef2}
A number q with trivial triangular remainder regarding  $m$ is defined as $T(q)=\frac{q^2+q}{2}\equiv T(t)=\frac{t^2+t}{2}{\pmod {m}}$ and $q=t$.
\subsubsection{Conclusion:}
All the numbers $m > 1$ have Trivial triangular remainder numbers. If $T(q)=\frac{q^2+q}{2}<m$ then q is number with a trivial triangular remainder number. 
\subsubsection{Example:}
The numbers with trivial triangular remainder regarding 15 are (0,1,2,3,4). 

\subsection{Definition of symmetric numbers with trivial triangular remainder}\label{Tdef3}
If $t$ is a number with trivial triangular remainder regarding $m$ then $s$ is a symmetric number with trivial triangular remainder if $s=m-t-1$. 
\subsubsection{Conclusion:}
There are the same quantity of symmetric numbers with  trivial triangular remainder than numbers with trivial triangular remainder.
\subsubsection{Example:}
The symmetric numbers with trivial quadratic number regarding 15 are (14,13,12,11).

\subsection{Definition of numbers with non trivial triangular remainder}\label{Tdef4}
A number q with non trivial triangular remainder regarding a number $m$ is defined as $T(q)=\frac{q^2+q}{2}\equiv T(d)=\frac{d^2+d}{2}{\pmod {m}}$ where $\frac{d^2+d}{2}<m$, $q\neq d$ and $m-q-1\neq d$
\subsubsection{Example:}
The numbers with non trivial triangular remainder regarding 15 are (5,6,8,9).

\section{Invariant and anti-invariant definition}
A complete description about invariant and anti-invariant definition can be found in \cite{Invariante}. Please refer to this article for more information.

\section{Properties of invariants}

\subsection{Proposition 1. Numbers and symmetric number with trivial and non trivial quadratic remainder can be inferred from the anti-invariant invariant tuples in odd numbers }\label{prop1}

\subsubsection{Demonstration:}
The tuples of anti invariant invariant number can be expressed as $(c,d)$ where $d=c+1$. If both numbers are multiplied  by $2\cdot s$, where $s^{2}<m$, then, are obtained two numbers  $c'=2\cdot c \cdot s$ and $d'=2\cdot c \cdot s+2\cdot s$. The equidistant number between c' and d' can be expressed as $e=2sc+s$ or also $e=2sd-s$. then $e^{2} {\pmod {m}}$ is:
\begin {equation}
e^2=4s^{2}cd+2ds^2-2cs^2-s^2\equiv 2ds^2-2cs^2-s^2{\pmod {m}}
\label{Peq1}
\end{equation}
if d is substituted by $d=c+1$ then (\ref{Peq1}) is now:
\begin {equation}
e^2\equiv 2cs^2+2s^2-2cs^2-s^2\equiv s^2{\pmod {m}}
\end{equation}


\subsection{Proposition 2. The trivial (0,1) anti-invariant invariant tuple inferred only numbers with trivial quadratic remainder in odd numbers}\label{prop2}
This is a particular case of \ref{prop1}. 
\subsubsection{Demonstration:}
The anti-invariant is 0. Then $e=2sc+s$ is always $e=s$. Due to $s^{2}<m$ there are not numbers with non trivial quadratic remainder because never meet $e\neq s$ . For this reason, a number with non trivial quadratic remainder cannot be inferred from (0,1) trivial anti-invariant, invariant tuple.

\subsection{Proposition 3. The trivial (m-1,m) anti-invariant invariant tuple inferred only symmetric numbers with trivial quadratic remainder in odd numbers}\label{prop3}
This is a particular case of \ref{prop1}. 
\subsubsection{Demonstration:}
 $(m-1,m) \equiv (-1,0) {\pmod {m}}$. The invariant is 0. Then $e=2sd-s$ is always $e=-s$ Due $s^{2}<m$, there are not number with  non trivial quadratic remainder because never meet $m-e\neq s$ . For this reason a number with non trivial quadratic remainder cannot be inferred from (m-1,m) trivial anti.invariant, invariant tuple.

\subsection{Proposition 4. The numbers with non trivial quadratic remainder are compounds numbers}\label{prop4}
The existence of a non trivial quadratic number probe automatically that the number is compound number. 
\subsubsection{Demonstration:}
Using \ref{prop2}, \ref{prop3} and the conclusions of \cite{Invariante} is probed that the numbers with non trivial quadratic remainder can be inferred only from the non trivial anti-invariant, invariant tuple. For \cite{Invariante} is demonstrated  that is a compound number.

\subsection{Proposition 5. The prime number and the remainder 1}\label{Nprop5}
If the number m is prime or power of prime $n^2\equiv 1$ only if n=1 or n=m-1.
\subsubsection{Demonstration:}
Is a particular case of \ref{prop4} where s=1. If exists a non trivial quadratic remainder 1, then exists a non trivial invariant and m is not prime.


\section{The triangular number and the invariant}

\subsection{Proposition 6. The invariant and the triangular number remainder.} \label{prop7}
The invariant also exists in the triangular numbers remainders. An Invariant A regarding m is a number where $T_{A}\equiv A{\pmod {m}}$.
\subsubsection{Demonstration:}
The multiplication of the tuple anti-invariant, invariant $(c,d)$ produces a rhomboid number. If this number is divide by 2 the result  is the triangular number, $T_{n-1}$ and multiple of m, therefore $T_{n-1}\equiv 0{\pmod {m}}$. the next triangular number:
\begin{equation}
T_{n}=T_{n-1}+n
\end{equation}
Then:
\begin{equation}
T_{n}=T_{n-1}+n\equiv n {\pmod {m}}
\end {equation}
\subsubsection{Conclusion:}
The invariant is also present in triangular number's remainder . There are tuples (d-1,d) where d-1 is zero remainder and d is an invariant. Some properties can be also extrapolate to the triangular number's remainder.

\subsection{Proposition 7. The trivial and non trivial triangular remainder can be inferred from the zero, invariant tuple in odd numbers }\label{prop8}

From the zero invariant tuple, the trivial and non trivial triangular remainder can be inferred. The number $n=2cs+c+s$ produces a triangle remainder if  $\frac{(s^2+s)}{2}<m$, where c is the position of the zero remainder.\\

\subsubsection{Demonstration:}

\begin{equation}\label{eqtri1}
T_{n}=\frac{(2cs+c+s)^2+2cs+c+s}{2}=\frac{4(c^2+c)(s^2+s)+(c^2+c)+(s^2+s)}{2}
\end{equation}

From \ref{prop7} the $\frac{c^2+c}{2}\equiv 0\;\;\;(mod\;m)$ then :\\
\begin{equation}
\frac{4(c^{2}+c)(s^{2}+s)+(c^{2}+c)+(s^{2}+s))}{2} \equiv \frac{(s^2+s)}{2} {\pmod {m}}
\end{equation}

\subsubsection{Conclusion:}
If $\frac{(s^2+s)}{2}<m$ all the numbers with trivial and non trivial triangular remainder can be inferred from the zero, invariant tuples.

\subsection{Proposition 8. The trivial (0,1) zero invariant tuple inferred only trivial triangular remainder numbers in odd numbers}\label{prop9}
This is a particular case of \ref{prop8}. 
\subsubsection{Demonstration:}
Is Supposed that in \ref{prop8} c=0, for this reason  n=s. Due to $\frac{(s^2+s)}{2}<m$ then n is a trivial triangular remainder. 

\subsection{Proposition 9. The trivial (m-1,m) zero invariant tuple inferred only trivial symmetric triangular remainder numbers in odd numbers}\label{prop10}
This is a particular case of \ref{prop8}. 
\subsubsection{Demonstration:}
$(m-1,m) \equiv (-1,0){\pmod {m}}$. For this c=-1 and  $n=-s-1$. The triangular remainder of  n is:
\begin{equation}
\frac{(-s-1)^2-s-1}{2}\equiv \frac{s^2+s}{2}{\pmod {m}}
\end {equation}
Then $n\equiv m-s-1{\pmod {m}}$ and $\frac{(s^2+s)}{2}<m$. As conclusion n is a trivial symmetric triangular remainder .

\subsection{Proposition 10. The numbers with non trivial triangular remainder are the compounds numbers}\label{prop11}
The existence of a non trivial triangular remainder probe automatically that the number is compound number 
\subsubsection{Demonstration:}
For \ref{prop9}, \ref{prop10} and the conclusions of \cite{Invariante} is probed that the non trivial triangular remainder can be inferred only from the non trivial zero, invariant tuple. The existence of a non trivial invariant (see \cite{Invariante}) demonstrate that the number is a compound number.

\section{Conclusions. The invariants properties can be used to develop new algorithms to test the primality }
From \ref{Nprop5}, \ref{prop11} and \cite{Invariante} is possible to play with different strategies in order to test the primality or get additional information to factorize the number.
The main strategies are:\\
\subsection{Using the invariant to test the primality}
The search for invariant could be used to test the primality and even provide additional information related to the factorization. The main advantage is that it is not necessary to check if the number is a power or triangular. On the other side is the method with the minimal potential solutions. There are $2^{\beta(m)}-2$ non-trivial invariants.\footnote{$\beta(m)$ is the number of prime factor of m.}\\
Additional to this, The algorithms can be distributed in different machines, can run in both directions and statistically strategies can be used.

\subsubsection{Algorithm 1. Algorithm to test the primality using invariants remains}
A very easy algorithm can be created to check the primality and even identifying if is a power of a prime or compound number just searching invariants. In order to do it is quite simple to create one algorithm. This algorithm go from down to up but could be done also in the opposite direction. The procedure is as following:
\begin{enumerate}
    \item Create one counter C1. This counter go from 2 to (m-1)/2 and each iteration is increased by 1.
    \item Create a second counter C2 this counter start in 4 and each iteration is increased by 2*C1-1.
    \item If C2 in one iteration is bigger than m then C2=C2-m.
    \item If C2=m then the program stop and m is a power.
    \item If C2=C1 then the program stop and m is a compound number.
    \item Finally if C1=(m-1)/2 an the program doesn´t stop before the program stop and is a prime number.
\end{enumerate}
As example is shown three different number and the result with the algorithm. The example is started with 55 \\

\begin{tabular}{| p{2cm}| p{2cm}| p{2cm}|p{4cm}|}
\hline
\textbf{Iteration}&\textbf{C1}&\textbf{C2}&\textbf{Result}\\ \hline
1	&2	&4 & \\ \hline
2	&3	&9 & \\ \hline
3	&4	&16&\\ \hline
4	&5	&25&\\ \hline
5	&6	&36&\\ \hline
6	&7	&49&\\ \hline
7	&8	&64-55=9&\\ \hline
8	&9	&26&\\ \hline
9	&10	&45&\\ \hline
10 &11 &66-55=11&\textbf{11=11 is compound}\\ \hline
\end{tabular}
\\
\\
If the same process is realized for a prime number like 23:\\

\begin{tabular}{| p{2cm}| p{2cm}| p{2cm}|p{4cm}|}
\hline
\textbf{Iteration}&\textbf{C1}&\textbf{C2}&\textbf{Result}\\ \hline
1	&2	&4 &\\ \hline
2	&3	&9 &\\ \hline
3	&4	&16 &\\ \hline
4	&5	&25-23=2 &\\ \hline
5	&6	&13 &\\ \hline
6	&7	&26-23=3&\\ \hline
7	&8	&18 &\\ \hline
8	&9	&35-23=12&\\ \hline
9	&10	&31-23=8 &\\ \hline
10 &11 &29-23=6 &\textbf{is prime}\\ \hline
\end{tabular}
\\
\\
Is observed that $C1=\frac{m-1}{2}$ and C2 is not 0 or equal to C1 then is prime
Finally if 9 is tested then:\\

\begin{tabular}{| p{2cm}| p{2cm}| p{2cm}|p{4cm}|}
\hline
\textbf{Iteration}&\textbf{C1}&\textbf{C2}&\textbf{Result}\\ \hline
1	&2	&4 &\\ \hline
2	&3	&9-9=0 &\textbf{At least one of his factor is a power.}\\ \hline
\end{tabular}
\\
\\
\textbf{Factorization}\\
The result provide a anti invariant invariant tuple (c,d) Each of then is multiple of one factor and for this reason factorizing c or d provide one of the factor numbers. Other strategy is multiply $c\cdot d=f$ then $f/m=g$ an the factorization of g give as result all the factors independents of c and d.\\ \\
\textbf{Code}\\
Following code implement the previous algorithm. This code is not optimize but give an idea about how easy is.

\begin{lstlisting}[language=Python]


 string = input('please insert a odd number:')
    num = int(string)
    # initialize control variable
    control = int((num - 1) / 2)
    prime = True
    # initialize the other variables
    C1 = 2
    C2 = 4
    # control loop
    while (control >= C1):

        if (C2 > num):
            C2 = C2 - num
        elif (C2 == num):
            print(f'the number is  raised to the a power{C1}:')
            prime = False
            break

        if (C1 == C2):
            print(f'the number is not a prime {C1}:')
            prime = False
            break
        C1 = C1 + 1
        C2 = C2 + 2 * C1 - 1
    if (prime):
        print(f'the number is a prime:')
\end{lstlisting}

\subsection{Using the non trivial quadratic remains to test the primality} 
Another possibility  could be the non-trivial quadratic remainders, which would also provide additional information related to the factorization. There are $\epsilon(m)(2^{\beta(m)}-2) $ solutions.\footnote{$\epsilon(m)$ is the quantity of numbers with trivial quadrangular remainder where GCD(n,m)=1.}\\ The main disadvantage is the process to check if it is a quadratic number.
These algorithms can be distributed on different machines and run in both directions. 
\subsubsection{Algorithm 2. Algorithm to test the primality using non trivial quadratic remains up down}
In this code the primality is tested searching non trivial quadratic remains in a number but from up to down. the algorithm is :
\begin{enumerate}
    \item Create one counter C1. this counter increase from  1 to (m-1)/2-$\sqrt{m}$ and each iteration is increased by 1.
    \item Create a second counter C2 this counter start in remains of the $((number-1)/2)^{2}$ and each iteration is increased by 2*C1.
    \item Create a number P2 where contains the square number nearest to C2. 
    \item If C2 in one iteration is bigger than m then C2=C2-m.
    \item If C2=m then the program stop and m is a power.
    \item If C2=P2 then the program stop and m is a compound number.
    \item Finally if C1=(m-1)/2-$\sqrt{m}$ an the program doesn´t stop before the program stop and is a prime number.
\end{enumerate}

For example the primality of 93 are going to be checked:
\\

\begin{tabular}{| p{3cm}| p{2cm}| p{2cm}|p{4cm}|}
\hline
\textbf{Iteration P1}&\textbf{C2}&\textbf{P2}&\textbf{Result}\\ \hline
0	&70	&81 &\\ \hline
1	&72	&81 &\\ \hline
2	&76	&81 &\\ \hline
3	&82	&100 &\\ \hline
4	&90	&100 &\\ \hline
5	&100-93=7	&100 &\\ \hline
6	&19	&25 &\\ \hline
7	&33	&36 &\\ \hline
8	&49	&49 &49=49 is not prime\\ \hline
\end{tabular}
\\
\\
If the same process is realized for a prime number like 23:\\

\begin{tabular}{| p{3cm}| p{2cm}| p{2cm}|p{4cm}|}
\hline
\textbf{Iteration P1}&\textbf{C2}&\textbf{P2}&\textbf{Result}\\ \hline
0	&6	&9 &\\ \hline
1	&8	&9 &\\ \hline
2	&12	&16 &\\ \hline
3	&18	&25&\\ \hline
4	&26-23=3&4 &\\ \hline
5	&13	&16&\\ \hline
6	&25-23=2	&18 &end is found is prime\\ \hline
\end{tabular}
\\
\\

Finally if 9 is tested then:\\

\begin{tabular}{| p{3cm}| p{2cm}| p{2cm}|p{4cm}|}
\hline
\textbf{iteration P1}&\textbf{C2}&\textbf{P2}&\textbf{result}\\ \hline
0	&7	&9 &\\ \hline
2	&9	&9-9=0 &\textbf{At least one of his factor is a power}\\ \hline
\end{tabular}
\\
\\
\textbf{Factorization}\\
This produce two square numbers $((m-1-2P1)/2)^{2}$ and P2. Obtaining the following numbers $c=((m-1-2P1)/2)+\sqrt{P2}$ and $d=((m-1-2P1)/2)-\sqrt{P2}$. For the properties of the invariant ca be seen that each contain a factor of m.

\subsection{Using the non trivial triangular remains to test the primality} 
Another possibility is to use non-trivial triangular remains, there are $\gamma(m)(2^{\beta(m)}-2)$ solutions.  \footnote{$\gamma(m)$ is the quantity of numbers with trivial triangular remainder where GCD(n,m)=1.}. This is the maximum number of possible solutions for all methods. The main disadvantage is the process of checking if it is a triangular number.

\subsubsection{Algorithm to test the primality using non trivial triangular remains}
In this code the primality is tested searching non trivial triangular remains in a number in direction from up to down. 
The algorithm is as following:
\begin{enumerate}
    \item Create one counter C1. This counter increase from  1 to (m-1)/2-s where s is the last number
    where $T(s)<m$ and each iteration is increased by 2.
    \item Create a second counter C2. This counter start in remains of the $(T(number-1)/2)$ and each iteration is increased by 2*C1.
    \item Create a third counter C3. This counter start in remains of the $(T(number-3)/2)$ and each iteration is increased by 2*(C1-1).
    \item Create a number P2 where contains the triangular number nearest to C2. 
    \item Create a number P3 where contains the triangular number nearest to C3.
    \item If C2 in one iteration is bigger than m then C2=C2-m.
    \item If C3 in one iteration is bigger than m then C3=C3-m.
    \item If C2=m then the program stop and m is a power.
    \item If C3=m then the program stop and m is a power.
    \item If C2=P2 then the program stop and m is a compound number.
    \item If C3=P3 then the program stop and m is a compound number.
    \item Finally if C1=(m-1)/2-s an the program doesn´t stop before the program stop and is a prime number.
\end{enumerate}

\begin{tabular}{| p{2cm}| p{2cm}| p{2cm}|p{2cm}|p{2cm}|p{2cm}|}
\hline
\textbf{iteration C1}&\textbf{C2}&\textbf{P2}&\textbf{C3}&\textbf{P3}&\textbf{result}\\ \hline
0	&12	&15 &58 &66 &\\ \hline
1	&16	&21 &60 &66 &\\ \hline
2	&24	&28 &66 &66 & T(45-2)= 946 and 66=66 is not prime \\ \hline
\end{tabular}
\\
\\
If the same process is realized for a prime number like 23:\\

\begin{tabular}{| p{2cm}| p{2cm}| p{2cm}|p{2cm}|p{2cm}|p{2cm}|}
\hline
\textbf{iteration C1}&\textbf{C2}&\textbf{P2}&\textbf{C3}&\textbf{P3}&\textbf{result}\\ \hline
0	&9	&10 &20 &21 &\\ \hline
1	&13	&15 &22 &25 &\\ \hline
2	&21	&21 &5 &6 & T(10-4)=21 21=21 is prime \\ \hline
\end{tabular}
\\
\\

Finally if 15 is tested then:\\

\begin{tabular}{| p{2cm}| p{2cm}| p{2cm}|p{2cm}|p{2cm}|p{2cm}|}
\hline
\textbf{iteration C1}&\textbf{C2}&\textbf{P2}&\textbf{C3}&\textbf{P3}&\textbf{result}\\ \hline
0	&6	&6 &13 &15 &T(6)=21 is not prime 6=6\\ \hline
1	&10	&10 &15 &15 & C3=number is a triangular number\\ \hline
2	&24	&28 &66 &66 & 66=66 is not prime \\ \hline
\end{tabular}
\\
\\

\end{document}